\documentclass{article}
\usepackage{comment}
\usepackage{commandesart}
\begin{document}
\title{The category of representations of a crossed module}
\date{}
\author{Ony Aubril\\}
\maketitle

\begin{abstract}
In this article, we present the notion of linear representations of an action groupoid induced by a group action underlying a crossed module of finite groups. We begin by defining action groupoids, then compute a few properties that are essentially tied to the action groupoid itself. We consider a functor category from this action groupoid to the category of finite dimensional vector spaces to link the notion to a previously known notion of representation of a crossed module of finite groups. From the latter, we inherit a braided tensor product. We prove that both notions are equivalent, then we show that the category of representations is actually a category of modules over a finite dimensional Hopf algebra. This construction leads to a new class of quasi-triangular Hopf algebras, indexed by crossed modules of finite groups. We prove that in the case of complex representations, this class of algebras is actually a class of semi-simple algebras, leading to a full description of simple subalgebras for any Hopf algebra of this class. Once this is established, we then study the monoidal structure of representations of a crossed module of finite groups, in particular links between direct sums of simple objects and tensor products of simple objects. We give results for a Clebsch-Gordan formula and in this context, detail proofs and develop the character theory for representations of crossed modules of finite groups. Lastly, we apply our results to the construction of ribbon invariants.
\end{abstract}

\section*{Introduction}

In various domains of mathematics, the concept of crossed modules naturally appears. Initially introduced in the 1940's by J.H.C. Whitehead in the context of groups \cite{Whitehead}, it originally was used as a tool to classify topological spaces up to homotopy via fundamental groups. A similar description appeared in the work of A.M. Turing in 1938 in a seemingly unrelated context \cite{Turing}. A.M.Turing introduced what would later be called a crossed module of groups as an efficient way to classify central group extensions given a kernel and cokernel of the extension, settling Hoelder \textit{extension problem} for groups. It is nowadays well-understood that crossed module of groups are a crucial tool to understand the categorical behaviour of the category of groups. Indeed, it has been shown that the category of crossed module of groups is equivalent to the category of internal categories in the category $\cgrp$ of groups and group morphisms (see section 8. of chapter XII in \cite{CategoriesFor}), or equivalently groups in the category of small categories (see \cite{MR0419643}), or even the category of $\text{Cat}^1\mhyphen$groups \cite{Spaceswithfinitelymanynon-trivialhomotopygroups}, showing that the category of crossed modules is actually a variety of algebras in the sense of F.W.Lawvere \cite{Functorialsemanticsofalgebraictheories}. This sequence of equivalences is referred to as the Brown-Spencer theorem. It was soon noticed that the notion of crossed module could be generalized to Lie algebras \cite{ExtensionsCentralesdAlgebresDeLie}, and then for cocommutative Hopf algebras \cite{Strictquantum2groups, Cat1Hopfalgebraandcrossedmodules}. Finally, the definition of crossed modules in any semi-abelian category in 2003, due to G. Janelidze, \cite{Janelidze}, was introduced. In this article, G. Janelidze proves that the Brown-Spencer theorem, stating that the category of internal groupoids is equivalent to the category of crossed modules of the same category, holds for any semi-abelian category, extending the theory of crossed modules to a large class of categories. The category of crossed modules for an algebraic structure is notably complex, and understanding it has led to various interpretations. The connections with the third cohomology group, for example, is recurrent within crossed modules of different algebraic structures (see \cite{WagXmodHn, W}), which gives an easier, and equivalent (see \cite{Clovis}), way to compute crossed module in the operadic context. Crossed modules of groups have also been shown to appear in the study of 2-categories, in a way similar to the way the group $\Aut{A}$ appears in a 1-category for each of its objects $A$ (see \cite{MonikaTruong}).

Our focus here is on the notion of crossed modules of finite groups. A theory of representations for such objects was first introduced by P. Bantay in \cite{B} as an example of premodular category. This approach notably includes a character theory that is a direct generalization of the character theory of group representations. This character theory has also been studied in \cite{DeghaniDavvaz}. Given that representation theory of groups is itself a complex subject, this generalization faces similar challenges, along with new ones. In this article we show that Bantay's category of representations $\mM(\xX)$, for a given crossed module of finite groups $\xX$, is actually a functor category from the action groupoid induced by $\xX$ to the category of finite dimensional vector spaces (proposition \ref{Equiv1-Cat}). This leads to a very efficient proof for the semi-simplicity of the category $\mM(\xX)$ over the complex field.

With this categorical background established, we then prove as one of the main results that $\mM(\xX)$ is equivalent as a monoidal category to a category of modules over some quasi-triangular Hopf algebra $\textsf{D}(\xX)$ that we fully describe in section \ref{Sect2}. As such, understanding the monoidal category $\mM(\xX)$ becomes equivalent to understanding the Hopf algebra $\textsf{D}(\xX)$. We give a formula for the number of simple $\textsf{D}(\xX)\mhyphen$modules up to isomorphism, and then a full description of simple modules up to isomorphism. Once this classification is established, we investigate relations between tensor products of simple objects and their decomposition into a direct sum of simple objects. We characterize simple objects that appear in the decomposition of some tensor product of simple objects in proposition \ref{Cleb-Gor}, and then consider the character theory that was given by P. Bantay in \cite{B}. We link his results to ones we established and develop further his theory by introducing a notion of character table very similar to the classical notion of character table for representations of finite groups, and which is a direct generalization of the latter. We prove that character tables must be diagonal by block and contain character tables of stabilizers of the group action in its diagonal (see \ref{char}). We give two examples of such character tables.

This new class of Hopf algebras opens up the search for ribbon invariants defined via representations of crossed modules of finite groups (\cite{Turaev}, chapter I, theorem 2.5, \cite{Ohtsuki} theorem 4.5. page 72). Ribbon invariants are especially interesting because it has been shown that they lead to invariants of 3-manifolds (\cite{Ohtsuki} chapter 8 for example), as 3-manifolds are obtained via surgery along ribbons. In the last section, we discuss some more properties of the category $\mM(\xX)$, or equivalently, the Hopf algebra $\textsf{D}(\xX)$. In view of finding quantum invariants of $3\mhyphen$manifolds, we prove that $\mM(\xX)$ is actually a ribbon category, and give explicit values for the ribbon invariant we get via Rechetikhin-Turaev\rq{}s theorem \cite{Turaev} in section \ref{LastSect}.

This work fully translates the foundations of finite dimensional representations of a given finite crossed module of groups into a theory of finite dimensional modules over an explicit Hopf algebra. This work also gives a detailed description of the behaviour of any Hopf algebra of the form $\textsf{D}(\xX)$ for some crossed module of finite groups $\xX$. This class of quasi-triangular Hopf algebras indexed by crossed modules of finite groups most notably includes the group algebra of any finite group, its dual algebra, as well as the Drinfeld quantum double of a finite group, among others.

The next step in the study of finite dimensional representations of crossed modules of finite groups would be to determine the extent of its applications. It is a good example of premodular category (see \cite{B}), of fusion, and of ribbon category. Its applications in quantum topology has yet to be studied more extensively, especially in conformal field theory where it may have a place as the generalization of the Drinfeld quantum double. 

\section*{Acknowledgments}
I would first like to greatly thank Friedrich Wagemann for leading this research project, his availability, and his guidance. I also am deeply thankful to Christian Blanchet for showing us that our category of interest was most likely a ribbon category, as well as its leads and contribution for ribbon invariants. I also wish to thank Christoph Schweigert for its feedback on the first version of this article, that lead to the consideration of action groupoids. 

\section{Action groupoids and Bantay\rq{}s representation of crossed modules}

We begin by recalling a few notions on $G\mhyphen$sets, $G\mhyphen$groups and crossed modules of finite groups.
\subsection{The category of $G\mhyphen$sets}

\begin{Definition}[$G\mhyphen$Sets]
If $G$ is a group, a $G\mhyphen$set is a pair $(X,\mu:G\times X\rightarrow X)$, where $X$ is a set and $\mu$ is a group action of $G$ on $X$. We define morphisms of $G\mhyphen$sets $f:(X,\mu)\rightarrow (Y,\nu)$ as set-theoretic maps $f:X\rightarrow Y$ such that $f(\mu(g,x))=\nu(g,f(x))$ for all $g\in G$, $x\in X$. The class of $G\mhyphen$sets together with the class of $G\mhyphen$sets morphisms define a category $\cgsets{G}$.
\end{Definition}

We shall usually write $g\cdot x$ for the action of $g\in G$ on $x$ in any $G\mhyphen$set, relying on context to interpret which action is involved.

The category $\cgsets{G}$ is actually quite well-behaved. It is an example of elementary topos (see \cite{SheavesinGeometryandLogic}). We can also think of it as the category of functors from $G$ considered as a one-object category to the category $\cset$ of sets and maps.
Throughout this article, we will mostly focus on $G\mhyphen$groups, that is, group objects in the category $\cgsets{G}$. The category $\cggrp{G}$ of $G\mhyphen$groups is the subcategory of $\cgsets{G}$ containing groups endowed with a $G$ action via group automorphisms, and morphisms are $G\mhyphen$set morphisms that also are group morphisms.

\begin{Definition}[Action groupoid]

Let $(X,\mu)$ be a $G\mhyphen$set. We define its associated action groupoid as the small category

\begin{center}
$X//G=$
\begin{tikzcd}
G\times G\times X \ar[r, "m_G\times \id{X}"]&G\times X\ar[r, shift left=4, "p_2"]\ar[r, shift right=4, "\mu"']&X\ar[l, "i_2" description]
\end{tikzcd}
\end{center}

where $m_G$ is the multiplication law of the group $G$. One can easily check that it is in fact a groupoid, with inverse morphism $s:G\times X\rightarrow G\times X$ given by $s(g,x)=(g\mm, g\cdot x)$. This assignment extends into a functor $F$ that sends a morphism of $G\mhyphen$sets $f:(X,\mu)\rightarrow (Y,\nu)$ to the functor of small categories

\begin{center}
\begin{tikzcd}
G\times G\times X \ar[r, "m_G\times \id{X}"]\ar[d, "\id{G\times G}\times f"']&G\times X\ar[d, "\id{G}\times f"']\ar[r, shift left=4, "p_2"]\ar[r, shift right=4, "\mu"' description]&X\ar[d, "f"]\ar[l, "i_2" description]\\
G\times G\times Y \ar[r, "m_G\times \id{Y}"]&G\times Y\ar[r, shift left=4, "p_2" description]\ar[r, shift right=4, "\nu"']&Y\ar[l, "i_2" description]
\end{tikzcd}
\end{center}
\end{Definition}
The notation $X//G$ is to understand as a weakening of $X/G$ the quotient of the set $X$ by the action of $G$. Indeed, instead of considering element of $X$ as equals when they are in the same orbit, we consider them as isomorphic. 

The functor $F$ is actually an embedding of categories from $\cgsets{G}$ to $\cCat$ the category of small categories. Indeed, $F$ is injective on objects (up to isomorphism) as $X//G\cong Y//G$ implies $X\cong Y$ in $\cset$ and this bijection commutes with the $G\mhyphen$action, so it is actually a morphism in $\cgsets{G}$, which makes it an isomorphism in $\cgsets{G}$. It is also faithful as $F(g)=F(f):X//G\rightarrow Y//G$ implies $F(f)$ is equal to $F(g)$ on objects and thus $f=g$.

If $X$ is a $G\mhyphen$group and not just a $G\mhyphen$set, we can define an action groupoid as before, which we will still write $X//G$. We now have all we need to state the first result, which will come back later on.

\begin{Proposition}\label{coprodStab}
If $X$ is a $G\mhyphen$set, then a skeleton of $X//G$ is the coproduct of categories $\displaystyle\coprod_{[x]\in X/G}\textbf{B}\Stab{G}{x}$, where $\textbf{B}\Stab{G}{x}$ is the one object category corresponding to the group $\Stab{G}{x}$.
\end{Proposition}
\begin{proof}
We already stated that $X//G$ is a groupoid. Considering this, a skeleton of $X//G$ is a disjoint union of one-object categories. Now if we look at hom-sets, $\Hom{X//G}{x,y}=\{(g,x)\vert g\cdot x=y \}$ is empty if $y\not \in G\cdot x$. If $y=h\cdot x$ for some $h\in G$, then $\Hom{X//G}{x,y}=\{(g,x)\vert x=(g\mm h)\cdot x \}\cong \Stab{G}{x}$. The equivalence of categories follows immediately.
\end{proof}
\subsection{First remarks on Bantay\rq{}s representations category}

We then get to the definition of our category of interest.

\begin{Definition}[Representation of $G\mhyphen$sets]
A representation of a $G\mhyphen$set $X$ (over a field $\Kk$) is a functor $X//G\rightarrow \cfd{\cvect{\Kk}}$ the category of finite dimensional vector spaces over a field $\Kk$. The category of representations of a $G\mhyphen$set $X$ is the functor category $\cFun{X//G}{\cfd{\cvect{\Kk}}}$.
\end{Definition}

For the rest of the article we will mostly be interested in $G\mhyphen$groups inducing a crossed module of finite groups. We recall the definition.
\begin{Definition}[Crossed module of groups]
A crossed module of groups is a quadruplet $\xX=(G,X,\mu, \partial)$ where $G$, $X$ are groups, $\mu$ is a left action of $G$ on $X$, written as $\mu(g,x)=g\cdot x$, via group automorphisms, and $\partial$ is a group morphism from $X$ to $G$ such that the following axioms are satisfied for all $x,y$ in $X$, $g$ in $G$ :
\begin{enumerate}
\item[\textbf{(PX)}] $\partial(g\cdot x)=g\partial(x)g\mm$
\item[\textbf{(Pei)}] $\partial(x)\cdot y=xyx\mm$
\end{enumerate}
\end{Definition}

In the literature, axiom \textbf{(PX)} is often called the precrossed module axiom, and \textbf{(Pei)} is referred to as the Peiffer axiom.

For example, for any group $G$ there is a crossed module $(G,G,\conj{\_},\id{G})$, where $\conj{\_}$ is the action by conjugation. Another crossed module that always exists when we have a group $G$ with neutral element $1$ is $(G, \{1\}, g\mapsto \id{\{1\}}, 1\mapsto 1 )$. This trivial example will turn out to be quite important as his representations will exactly be representations of the group $G$, embedding the theory of finite group representations in our theory of representations of a crossed module of finite groups. More generally, a group $G$ with normal subgroup $N$ also gives rise to a crossed module of groups $(G,N,\conj{\_}, i)$, where $i$ is the inclusion and $\conj{\_}$ is the conjugation by elements of $G$ on $N$. Lastly, to give a different kind of example, we also have, for any group $G$, a crossed module of groups $(\Aut{G}, G, \id{\Aut{G}}, g\mapsto \conj{g})$.

From a crossed module of finite groups $\xX=(G,X,\mu,\partial)$, we get an underlying $G\mhyphen$group $X$, with the action $\mu$. We say that the category of representations of a crossed module of finite groups is the category of representations of this underlying $G\mhyphen$group. In \cite{B}, P. Bantay gave another definition for the notion of representations of a crossed module of finite groups. We will show that both definitions are actually equivalent, at least as categories. 

From a crossed module of finite groups $\xX=(G, X, \mu, \partial)$ arises a category of representations over a field $\Kk$, introduced in \cite{B}, that we will write $\mM(\xX)$. Objects in this category are defined as follows.

\begin{Definition}[Bantay\rq{}s representation of a finite crossed module of groups, \cite{B}]
Let $\xX=(G,X\mu,\partial)$ be a crossed module of finite groups. A finite dimensional representation of $\xX$ is a finite dimensional $\Kk\mhyphen$vector space $V=\bigoplus_{x\in X} V_x$, together with projectors $P(x)$, such that $V$ is a representation of the finite group $G$. The action of $G$ is written as a group morphism $Q:G\rightarrow GL(V)$, with $GL(V)$ the group of automorphisms of the vector space $V$, and for all $x$ in $X$, for all $g$ in $G$, $$P(g\cdot x)=Q(g)P(x)Q(g\mm).$$
\end{Definition}

A morphism in $\mM(\xX)$ is a linear map that is required to preserve both the $X\mhyphen$grading and the action of $G$. One can easily see that this definition makes sense for any $G\mhyphen$set $X$. P. Bantay defined the category above for the tensor product he was able to get on it, which requires the structure of group over $X$, and for the braiding on this tensor product, which is tied to the morphism $\partial$ above. Before diving into these considerations, we give the following result that is about the category $\mM(\xX)$ itself, and does not account for the tensor product or braiding.

\begin{Proposition}\label{Equiv1-Cat}
There is an equivalence of categories $$\cFun{X//G}{\cfd{\cvect{\Kk}}}\cong \mM(\xX).$$
\end{Proposition}
\begin{proof}
To a functor $\gym$ in $\cFun{X//G}{\cfd{\cvect{\Kk}}}$ we associate the representation in the sense of Bantay whose underlying vector space is $\bigoplus_{x\in X}\gym(x)$, with projections $\{P(x)\vert x\in X\}=\{\gym(1,x)\vert x\in X\}$. The $G\mhyphen$action $Q$ is given by $(Q(g): \gym(x)\rightarrow \gym(g\cdot x))=\gym(g,x)$, which can be considered as a family of automorphisms of $\gym(X)$ by universal property of the biproduct.
We can then check that $Q(g)P(x)Q(g\mm)$ is non trivial if and only if it acts on the component $\gym(g\cdot x)$, and then $Q(g)P(x)Q(g\mm)(\gym(g\cdot x))=\gym(((g,x)\circ (1,x)\circ (g\mm, g\cdot x))(g\cdot x))=\gym((1, g\cdot x)(g\cdot x))=P(g\cdot x)(\gym(g\cdot x)).$ A natural transformation $\alpha$ between
a functor $\gym$ and another functor $\fym$ is sent to a family of linear morphisms $\{\alpha_x :\gym(x)\rightarrow \fym(x)\vert x\in X\}$, which is equivalent to a morphism of graded vector space $\alpha : \bigoplus_{x\in X}\gym(x)\rightarrow \bigoplus_{x\in X}\fym(x)$ by universal property of the biproduct.
Naturality of $\alpha$ is exactly the compatibility requirement with the $G\mhyphen$action.
\end{proof}

\begin{Corollary} \label{SemiSimple}
If we suppose $\Kk=\Cc$ and $\xX=(G,X,\mu,\partial)$ is a crossed module of finite groups, then $\mM(\xX)$ is a semi-simple category, with finite number of simple objects up to isomorphism.
\end{Corollary}
\begin{proof}
Using both propositions \ref{coprodStab} and \ref{Equiv1-Cat}, we get the equivalence of categories
\begin{align*}
\mM(\xX)&\cong \cFun{X//G}{\cfd{\cvect{\Cc}}}\\
&\cong \cFun{\coprod_{[x]\in X/G}\textbf{B}\Stab{G}{x}}{\cfd{\cvect{\Cc}}}\\
&\cong \prod_{[x]\in X/G}\cFun{\textbf{B}\Stab{G}{x}}{\cfd{\cvect{\Cc}}}.
\end{align*}
The last equivalence makes clear our category of interest is a finite product of categories of finite dimensional complex representations of finite groups, which are semi-simple by the classical Maschke’s theorem. The result is then immediate.
\end{proof}

In \cite{B}, along with the category $\mM(\xX)$, for a crossed module of finite groups $\xX$, P. Bantay introduced a tensor product $\ot$ defined as follows : for $V,W$ objects in $\mM(\xX)$,

$$\adap{\bigoplus_{x\in X}V_x, Q_V}\ot \adap{\bigoplus_{x\in X} W_x, Q_W}=\adap{\bigoplus_{x\in X}\adac{\bigoplus_{y\in X}V_{y}\ot_\Kk W_{y\mm x}},Q_V\ot_\Kk Q_W}.$$ The natural transformation whose $(V,W)$ component is given by
\begin{align*}
\mapp{\beta_{V,W}}{V\ot W}{W\ot V}{v\ot w}{\dsum_{x\in X}Q_W(\partial(x))w\ot P_V(x)v}
\end{align*}
is a braiding for $(\mM(\xX), \ot)$.

Now our previous interpretation fails to account for the structure of monoid over the category itself. A group by itself does not have the comultiplication law necessary to encode the above grading of the tensor product. We then turn to Hopf algebras in the next section to give an equivalence of monoidal categories.

\section{A Hopf algebra to represent a crossed module of groups}\label{Sect2}
\subsection{The construction}

In \cite{K}, C. Kassel explains how to form Hopf algebras from a pair of Hopf algebras and suitable actions. The Drinfeld quantum double of any Hopf algebra is such an example of Hopf algebra. In this section, we use some results from \cite{K} to construct an explicit quasi-triangular Hopf algebra $\textsf{D}(\xX)$ from a crossed module of finite groups $\xX=(G,X\mu, \partial)$. This Hopf algebra $\textsf{D}(\xX)$ will contain both $\Kk G$, the group algebra of $G$ over $\Kk$, and $\Kk^X$, the dual of $\Kk X$, as well as the action of $G$ on $X$. The construction is based on chapter IX of \cite{K}.

Let $\xX=(G,X\mu, \partial)$ be a crossed module of finite groups. The linearisation of $\mu$ turns $\Kk^X$ into a left $\Kk G$-module-algebra by defining the module structure as $a\cdot \kr{x}:=\kr{a\cdot x}$ for all $a$ in $G$ and $x$ in $X$. Here and throughout the article $\kr{x}$ will represent the dual of the element $x$ of $\Kk X$. Note that we also use $\kr{x,y}$ for the Kronecker symbol, that is, $\kr{x,y}=0$ if $x\not= y$ and $1$ otherwise. 

\begin{Proposition}\label{DefDX}
Let $\xX=(G,X,\mu, \partial)$ a crossed module of finite groups. $\textsf{D}(\xX):=(\Kk^X\ot \Kk G, m, \eta, \Delta, \varepsilon, S)$ is a Hopf algebra, with the following mappings given on the bases $G$ of $\Kk G$ and $\{\kr{x}\}_{x\in X}$ of $\Kk^X$ :
\begin{itemize}
\item The multiplication is given by $m(\kr{x}\ot a\ot \kr{y}\ot b)=\kr{x,a\cdot y}\kr{x}\ot ab$,
\item the unit is $\eta(1)=\sum_{x\in X}\kr{x}\ot 1$
\item the comultiplication is $\Delta(\kr{x}\ot a)=\sum_{y\in X}\kr{y}\ot a\ot \kr{y\mm x}\ot a$,
\item the counit is $\varepsilon(\kr{x}\ot a)=\kr{x,1}$,
\item lastly, the antipode is $S(\delta_x\ot a)=\delta_{a^{-1}\cdot x\mm}\ot a^{-1}$.
\end{itemize}
\end{Proposition}

\begin{Remark}\label{Involutive}
$\textsf{D}(\xX)$ is an involutive Hopf algebra, i.e. $S\circ S=\id{\textsf{D}(\xX)}$. For $\delta_x\ot a$ in $\textsf{D}(\xX)$, $S\circ S(\delta_x\ot a)=S(\delta_{a^{-1}\cdot x\mm}\ot a^{-1})=\delta_{a\cdot(a^{-1}\cdot x\mm)\mm}\ot a=\delta_x\ot a$.
\end{Remark}

We can now give the main theorem here, whose proof is straightforward.
\begin{theorem}\label{ModCat}
Let $\xX=(G,X,\mu,\partial)$ be a finite crossed module of groups. Then $\mM(\xX)$ and $\cfd{\cmod{\textsf{D}(\xX)}}$ are isomorphic as monoidal categories.
\end{theorem}
\begin{proof}
    As mere categories, the isomorphism comes from the following. Since $(\cfd{\cvect{\Kk}},\ot,\Kk)$ is a closed monoidal category, a $\textsf{D}(\xX)\mhyphen$module structure over a finite dimensional vector space $V$ can be thought of as a linear map $\phi:\textsf{D}(\xX)\to \Hom{\cfd{\cvect{\Kk}}}{V,V}$. If we write, for $x\in X$, $P_\phi(x):=\phi(\kr{x}\ot 1)$ and $Q_\phi(a):=\phi(\sum_{x\in X}\kr{x}\ot a)$, we get 
    \begin{align*}
        \id{V}&=\phi(\eta(1))=\phi(\sum_{x\in X}\kr{x}\ot 1)\\
        &=\sum_{x\in X}\phi(\kr{x}\ot 1)=\sum_{x\in X}P_\phi(x).
    \end{align*}
    Similarly, $P_\phi(x)P_\phi(y)=\phi((\kr{x}\ot 1)(\kr{y}\ot 1))=\phi(\kr{x,y}\kr{x}\ot 1)=\kr{x,y}P_V(x)$. Moreover, for all $a\in G$, the assignment $a\in G\mapsto Q_\phi(a)$ is a group morphism $G\to GL(V)$. It is clear that $Q_\phi(1)=\phi(\eta(1))=\id{V}$, and 
    \begin{align*}
        Q_\phi(a)Q_\phi(b)&=\phi((\sum_{x\in X}\kr{x}\ot a)(\sum_{y\in X}\kr{y}\ot b))\\
        &=\phi(\sum_{x,y\in X}\kr{x,a\cdot y}\kr{x}\ot ab)\\
        &=\phi(\sum_{x\in X}\kr{x}\ot ab)\\
        &=Q_\phi(ab).
    \end{align*}
    Lastly, 
    \begin{align*}
        Q_\phi(a)P_\phi(x)Q_\phi(a\mm)&=\phi(\sum_{z,y\in X}(\kr{z}\ot a)(\kr{x}\ot 1)(\kr{y}\ot a\mm))\\&=\sum_{y\in X}\phi(\kr{y,a\cdot x}\kr{y}\ot aa\mm)\\&=\phi(\kr{a\cdot x}\ot 1)=P_\phi(a\cdot x),
    \end{align*}
    and we built an object $(V,P_\phi,Q_\phi)$ in $\mM(\xX)$ from an object $(V,\phi)$ in $\cfd{\cmod{\textsf{D}(\xX)}}$. 
    Conversely, an object $(V,P,Q)$ in $\mM(\xX)$ induces a linear map $\textsf{D}(\xX)\to \Hom{\cfd{\cvect{\Kk}}}{V,V}$ by giving over the basis $\{\kr{x}\ot a\}_{x\in X, a\in G}$, $\phi_{P,Q}(\kr{x}\ot a):=P(x)Q(a)$. By basic properties of $P$ and $Q$ we can recover $P_{\phi_{P,Q}}=P$ and $Q_{\phi_{P,Q}}=Q$. Similarly $\phi_{P_\phi,Q_\phi}(\kr{x}\ot a)=P_\phi(x)Q_\phi(a)=\phi(\kr{x}\ot a)$, which proves that there is a bijection between objects in $\mM(\xX)$ and objects in $\cfd{\cmod{\textsf{D}(\xX)}}$. It is immediate that a linear map between two objects in $\mM(\xX)$ is compatible with operations $P$ and $Q$ if and only if it is compatible with the corresponding action $\phi_{P,Q}$. This observation turns the bijection into an isomorphism of categories. 

    The next step is to prove that this isomorphism respects the monoidal structure. Over the vector space itself, in $\mM(\xX)$ and in $\cfd{\cmod{\textsf{D}(\xX)}}$ the tensor product is nothing more than the usual tensor product of vector space, with unit $\Kk$ and usual associativity constraints. For two objects $(V,\phi), (W,\psi)$ in $\cfd{\cmod{\textsf{D}(\xX)}}$, for $x\in X$, 
    \begin{align*}
        P_{\phi\ot \psi}(x)&=(\phi\ot \psi)\Delta(\kr{x}\ot 1)\\&=(\phi\ot \psi)(\sum_{y\in X}\kr{y}\ot 1\ot \kr{y\mm x} \ot 1)\\&=\sum_{y\in X}P_\phi(y)\ot P_\psi(y\mm x),
    \end{align*} 
    and we recover the $X\mhyphen$graduation of the tensor product in $\mM(\xX)$. Similarly for $a\in G$, \begin{align*}
        (Q_{\phi\ot\psi})(a)&=(\phi\ot \psi)\Delta(\sum_{x\in X}\kr{x}\ot a)\\&=\sum_{x,y\in X}\phi(\kr{x}\ot a)\ot \psi(\kr{y}\ot a)\\&=Q_\phi(a)\ot Q_\psi(a),
    \end{align*} and the isomorphism of categories extends into an isomorphism $\cfd{\cmod{\textsf{D}(\xX)}}\cong \mM(\xX)$ of monoidal categories. 
\end{proof}
One can easily check that taking $X=\{1\}\subseteq G$ we find $\textsf{D}(\xX)\cong \Kk G$ the algebra of the group $G$. This remark justifies the inclusion of the theory of representations of finite groups into Bantay's theory of representations of crossed modules of finite groups.

The category of Bantay's representations of a crossed module of finite groups is known to be braided, therefore $\textsf{D}(\xX)$ becomes a quasi-triangular Hopf algebra. We get the $R\mhyphen$matrix via the following result.

\begin{Proposition}[\cite{K}]
Let $(A,m,\eta,\Delta,\varepsilon)$ be a bialgebra. Then $\cmod{A}$ is braided if and only if the bialgebra $A$ is quasi-triangular.
Moreover, if $c$ is the braiding in $\cmod{A}$, and $\tau$ is the flip, then the $R\mhyphen$matrix on $A$ is given by $R=\tau_{A,A}c_{A,A}(\eta(1)\ot \eta(1))$, with $A$ seen as a module over itself.
\end{Proposition}
\begin{Corollary}
$\textsf{D}(\xX)$ is braided, with $R$-matrix given by $R=\dsum_{x\in X}\delta_x\ot 1\ot 1\ot\partial(x).$
\end{Corollary}

\subsection{Simple subalgebras and the algebra structure of $\textsf{D}(\xX)$}

\subsubsection{Semi-simplicity of $\textsf{D}(\xX)$}

\begin{Proposition}
If the characteristic of $\Kk$ does not divide the order of $G$, then $\textsf{D}(\xX)$ as defined in \ref{DefDX} is a semi-simple associative algebra.
\end{Proposition}
\begin{proof}
This is a special case of theorem 10.4.3. in \cite{Rad}. An immediate corollary of this theorem is that every finite dimensional involutive Hopf algebra $A$ is semi-simple as an associative algebra if the characteristic of $\Kk$ does not divide the dimension of $A$ as a vector space. We can also consider corollary \ref{SemiSimple}, in light of theorem \ref{ModCat}. We give a third argument by exhibiting a projection onto an invariant complementary subspace. If $M$ is a $\textsf{D}(\xX)\mhyphen$module and $N$ is a sub$\mhyphen \textsf{D}(\xX)\mhyphen$module of $M$, let $p$ be a linear projection from $M$ to the sub-vector space $N$ of $M.$ The projection $P_N$ onto $N$ from $M$ as $\textsf{D}(\xX)\mhyphen$modules is given by
\begin{align*}
\mapp{P_N}{M}{N}{m}{\dfrac{1}{\vert G\vert }\dsum_{g\in G, h\in X}}(\kr{h}\ot g)\cdot p((\kr{g\mm\cdot h }\ot g\mm)\cdot m).
\end{align*}
One can prove this is indeed a fitting projection by using an integral on the Hopf algebra. Exhibiting this projection gives a slightly better result about the characteristic of the field $\Kk$ for $\textsf{D}(\xX)$ to be semi-simple.
\end{proof}
Notice that if $X$ is the trivial group, then we have the usual Maschke projection for group representations.
From now on we assume $\Kk=\Cc$, so $\textsf{D}(\xX)$ is semi-simple.

\subsubsection{The centre of $\textsf{D}(\xX)$ and formulas for simple subalgebras}

Direct computations lead to an expression for a basis of the centre of $\textsf{D}(\xX)$.
\begin{Proposition}
For a crossed module of finite group $\xX=(G,X,\mu,\partial)$, the family
$$\bB= \adaa{\dsum_{g \in G}\delta_{g\cdot z}\ot gcg^{-1}\in \textsf{D}(\xX) \;\mid\; z\in X, c\in \Stab{G}{z}}$$ is a basis of $Z(\textsf{D}(\xX))$ the centre of the algebra $\textsf{D}(\xX)$.
\end{Proposition}

Of course the indexing above is not optimal, some elements in $\bB$ can correspond to two different indexes. Next corollary gives a proper indexing for $\bB$. 

\begin{Corollary}
Let $\hH$ be a set of representatives of $X/G$, the set of orbits of $X$ under the action of $G.$
$$\di{Z(\textsf{D}(\xX))}=\card{\{(x, [a]) \vert x\in \hH, [a] \text{ a conjugation class in }\Stab{G}{x}\}}$$
\end{Corollary}
\begin{proof}
Let $\aA=\{(x, [a]) \vert x\in \hH, [a] \text{ a class of conjugation in }\Stab{G}{x}\}.$
We are going to exhibit a bijection between $\aA$ and $\bB.$

Consider the map $$\mapp{\psi}{\aA}{\bB}{(x, [a])}{\dsum_{g\in G}\delta_{g\cdot x}\ot gag^{-1}}.$$

First, we need to check that this is well defined.
Let $x\in \hH, a,b\in\Stab{G}{x}$ such that $[a]=[b].$ Then by definition there exists $c\in \Stab{G}{x}$ such that $a=cbc\mm $. Therefore, we have the following equalities,
\begin{align*}\dsum_{g\in G}\delta_{g\cdot x}\ot gag\mm =\dsum_{g\in G}\delta_{g\cdot x}\ot gc b (gc )\mm=\dsum_{g\in G}\delta_{(gc\mm )\cdot x}\ot gbg\mm=\dsum_{g\in G}\delta_{g\cdot x}\ot gbg\mm,\end{align*}
so $\psi$ is well defined.
Let us prove that it is surjective.

Let $\beta=\dsum_{g\in G}\delta_{g\cdot x}\ot gag\mm\in \bB.$ Let $\chi$ be the representative of $G\cdot x$ that is in $\hH.$ We can then write $x=k\cdot \chi$, for some $k\in G.$ We just have to check that $k\mm ak\in \Stab{G}{\chi}$, which is true because $(k\mm ak)\cdot \chi=k\mm\cdot x=\chi.$ Follows

\begin{align*}\psi(\chi,[k^{-1}ak])=\dsum_{g\in G}\delta_{g\cdot \chi}\ot gk\mm ak g\mm =\dsum_{g\in G}\delta_{(gk)\cdot \chi}\ot gag\mm=\beta.\end{align*}

Lastly, we prove injectivity.
Let $(x,[a]), (y,[b])\in \aA$ such that $\psi(x,[a])=\psi(y,[b])$, i.e.
\begin{align*}
\dsum_{g\in G}\delta_{g\cdot x}\ot gag\mm=\dsum_{g\in G}\delta_{g\cdot y}\ot gbg\mm.
\end{align*}
Every term of each sum is an element of the canonical basis of $\textsf{D}(\xX)$ as a vector space, a consequence is that for this equality to be true, each term of the sum needs to be exactly equal to one term of the other sum. In particular, $y$ needs to be in the orbit of $x$. Both being representatives, we have $x=y$. From this assertion, we can say that both $[a]$ and $[b]$ are conjugation classes in $\Stab{G}{x}.$ Rewriting the equality, we have
\begin{align*}\dsum_{g\in G}\delta_{g\cdot x}\ot g(a-b)g\mm=0\ri \dsum_{g\in \Stab{G}{x}}(gag\mm-gbg\mm)=0\ri [a]=[b].\end{align*}
\end{proof}

As an interesting consequence, we can give a proof for the cardinality of simple objects up to isomorphism in $\mM(\xX)$. Indeed, one can consider Artin-Wedderburn theorem (see for instance \cite{P}, theorem 3.5 page 49). As an associative algebra, $\textsf{D}(\xX)$ is isomorphic to a product of algebras of matrices with complex coefficients. Every such algebra is a simple algebra, and for every such algebra, the centre is of dimension $1$. In consequence, the dimension of the centre of an associative semi-simple algebra is exactly the number of simple sub-algebras, which is exactly the number of simple objects up to isomorphism in the category of modules. This following theorem already appeared in \cite{L}.

\begin{theorem}\label{nbrsimples}
Let $\xX=(G,X,\mu,\partial)$ be a crossed module of finite groups, and $\hH$ be a set of representatives for $X/G$. There is a one-to-one correspondence
$$
\adaa{\text{simple object in $\mM(\xX)$}}\overset{1-1}{\longleftrightarrow}\adaa{(h, \chi_h)\vert h\in \hH, \chi_h\text{ is an irreducible character of } \Stab{G}{h}}
$$
\end{theorem}

To be able to have an explicit form of the simple subalgebras of $\textsf{D}(\xX)$, we can exhibit the family of primitive central sum-1 orthogonal idempotents of $\textsf{D}(\xX)$. An isomorphism then comes from \cite{P}, proposition a., p.94. In the next subsection, we will use this result to deduce a Clebsch-Gordan formula (Proposition \ref{Cleb-Gor}) in our context. We obtain it again via character theory in section \ref{TheClebschGordanFormulaViaCharacterTheory}.

\begin{Proposition}\label{idemp}
Let $\xX=(G,X,\mu,\partial)$ be a crossed module of finite groups and $\hH$ a family of representatives for $G/X$. For each $z\in \hH,$ let $\{e_i^z\}_{1\leq i\leq n_z}$ be the family of central sum-1 primitive orthogonal idempotents for the group algebra $\Cc\Stab{G}{z}$. Then the family $$\fF=\adaa{f_i^z=\dfrac{1}{\vert \Stab{G}{z}\vert }\dsum_{g\in G}\kr{g\cdot z}\ot g e_i^zg\mm}_{z\in \hH, 1\leq i\leq n_z}$$ is a family of central primitive orthogonal sum-1 idempotents for the algebra $\textsf{D}(\xX).$
\end{Proposition}

\subsubsection{A Clebsch-Gordan formula for representations of a crossed module.}

We now know what simple representations of a given crossed module of finite groups $\xX=(G,X,\mu, \partial)$ are. We saw in the preliminaries that there exists a notion of tensor product in the category $\mM(\xX)$. We will now want to understand the decomposition into simple objects of the tensor product of two simple objects.

Simple modules are known thanks to the family $\fF$ above, and to see if a simple finite dimensional $\textsf{D}(\xX)\mhyphen$module appears in the decomposition of some other $\textsf{D}(\xX)\mhyphen$module $M$, we just have to see if the corresponding element of $\fF$ is "orthogonal" to $M$. In the following sections $I^a_i$ will be the simple subalgebra corresponding to the idempotent element $f_i^a$ in $\fF$, for $a\in \hH $ and $i\in \llbracket 1, n_a\rrbracket.$ Let us write $e_i^{g\cdot z}:=ge_i^zg\mm$. It is easy to see that the family $\{e_i^{g\cdot z}\vert 1\leq i\leq n_z\}$ is actually a family of central primitive orthogonal sum$\mhyphen 1$ idempotents elements of $\Cc\Stab{G}{g\cdot z}$.

Direct computations lead to the following proposition.
\begin{lemma}
If $I_k^t$ is a simple module that appears in the decomposition of $I_j^z\ot I_i^s$ in simple $\textsf{D}(\xX)\mhyphen$modules, then $t$ is an element of $(G\cdot z)(G\cdot s)$.
\end{lemma}
\begin{proof}
We can directly compute $(f_k^t)\cdot (f_j^z\ot f_i^s)$ and show that this element is zero if $t$ is not an element of $(G\cdot z)(G\cdot s)$. We write $e_k^t=\sum_{a\in \Stab{G}{t}}\lambda_aa$.
\begin{align}
\vert \Stab{G}{z}\vert& \vert \Stab{G}{s}\vert\vert \Stab{G}{t}\vert(f_k^t)\cdot (f_j^z\ot f_i^s)=(\dsum_{g\in G}\kr{g\cdot t}\ot e_k^{g\cdot t})\cdot (\dsum_{h,l\in G}\kr{h\cdot z}\ot e_i^{h\cdot z}\ot \kr{l\cdot s}\ot e_j^{l\cdot s})\nonumber
\\&=\sum_{g,h,l\in G, m\in X,a\in \Stab{G}{t}}\lambda_a(\kr{m}\ot gag\mm)(\kr{h\cdot z}\ot e_i^{h\cdot z})\ot (\kr{m\mm (g\cdot t)}\ot gag\mm)(\kr{l\cdot s}\ot e_j^{l\cdot s})\label{calcul}
\end{align}
The first multiplication in the sum is non-zero only if $m=(gag\mm h)\cdot z$, and looking at the second multiplication with this value of $m$, we get the Kronecker symbol $\kr{((ag\mm h)\cdot z\mm) t, (ag\mm l)\cdot s}$, which is zero if $t$ is not in $(G\cdot z)(G\cdot s)$.
\end{proof}

With this established, actually computing simple subalgebras can be done by only looking at the group algebra $\Cc G$.

\begin{Proposition}\label{Cleb-Gor}
For any $s,z,t$ in $X$, let $\hH_{z,s}^t=\{(h,l)\in G^2 \vert (h\cdot z)(l\cdot s)=t\}.$ Let us fix $s,z,t\in \hH$ a set of representatives for $X/G$. Let $(i,j,k)\in \llbracket1, n_z\rrbracket\times\llbracket1,n_s\rrbracket\times\llbracket1, n_t\rrbracket$.
Then $I_k^t$ is a sub$\mhyphen \textsf{D}(\xX)\mhyphen$module of $I_j^z\ot I_i^s$ if and only if both of the following are true.
\begin{itemize}
\item $t\in (G\cdot z)(G\cdot s).$
\item As sub$\mhyphen \Cc G\mhyphen$modules of $\Cc G$, $\bigoplus_{g\in G}e_k^{g\cdot t}\Cc \Stab{G}{g\cdot t}$ is a right-submodule of $$\bigoplus_{(h,l)\in \hH_{z,s}^{g\cdot t}} e_j^{h\cdot z}\Cc \Stab{G}{h\cdot z}\ot e_i^{l\cdot s}\Cc \Stab{G}{l\cdot s}.$$
\end{itemize}
\end{Proposition}
\begin{proof}
We go back to equation (\ref{calcul}). We get with a change of variable $(h_0, l_0)=(g\mm h, g\mm l)$, and using the Sweedler notation: $x'\ot x''=\Delta(x)$. 
\begin{align*}
(\text{\ref{calcul}})&=\sum_{g\in G, (h_0,l_0)\in \hH_{z,s}^t, a \in \Stab{G}{t}}\lambda_a\kr{(gah_0)\cdot z}\ot ga e_i^{h_0\cdot z}g\mm\ot \kr{(gah_0)\cdot z\mm (g\cdot t)}\ot ga e_j^{l_0\cdot s}g\mm \\
&=\sum_{g\in G, (h_0,l_0)\in \hH_{z,s}^t, a \in \Stab{G}{t}}\lambda_a\kr{(gah_0)\cdot z}\ot ga e_i^{h_0\cdot z}g\mm\ot \kr{(gal_0)\cdot s}\ot ga e_j^{l_0\cdot s}g\mm\\
&=\sum_{g\in G, (h,l)\in \hH_{z,s}^t, a \in \Stab{G}{t}}\lambda_a\kr{(gh)\cdot z}\ot g e_i^{h\cdot z}ag\mm\ot \kr{(gl)\cdot s}\ot g e_j^{l\cdot s}ag\mm\\
&=\sum_{g\in G, (h,l)\in \hH_{z,s}^t}(\kr{(gh)\cdot z}\ot e_i^{(gh)\cdot z}\ot \kr{(gl)\cdot s}\ot e_j^{(gl)\cdot s})(1\ot (e_k^{g\cdot t})\rq{}\ot 1\ot (e_k^{g\cdot t})\rq{}\rq{})\\
&=\sum_{g\in G, (h,l)\in \hH_{z,s}^{g\cdot t}}(\kr{h\cdot z}\ot e_i^{h\cdot z}\ot \kr{l\cdot s}\ot e_j^{l\cdot s})(1\ot (e_k^{g\cdot t})\rq{}\ot 1\ot (e_k^{g\cdot t})\rq{}\rq{})\\
&=\vert \Stab{G}{t}\vert\vert \Stab{G}{z}\vert \Stab{G}{s}\vert \sum_{\tau\in G\cdot t, \zeta\in G\cdot z, \sigma\in G\cdot s, \zeta\sigma=\tau}(\kr{\zeta}\ot e_i^{\zeta}\ot \kr{\sigma}\ot e_j^{\sigma})(1\ot (e_k^{\tau})\rq{}\ot 1\ot (e_k^{\tau})\rq{}\rq{})
\end{align*}

This last sum is a free sum for both $\zeta$ and $\sigma$, thus the entire sum is zero if and only if each term in $(\zeta,\sigma)$ is zero, which is equivalent as $\sum_{\tau\in G\cdot t}(e_i^\zeta\ot e_j^{\sigma})\Delta(e_k^\tau)=0$ for all couples $(\zeta, \sigma)$ in the sum. This is the case if and only if $\bigoplus_{g\in G}e_k^{g\cdot t}\Cc \Stab{G}{g\cdot t}$ is not a right submodule of $ e_j^{h\cdot z}\Cc \Stab{G}{h\cdot z}\ot e_i^{l\cdot s}\Cc \Stab{G}{l\cdot s}$ for any $(h,l)\in \hH_{z,s}^{g\cdot t}$. This concludes the proof.
\end{proof}

With this approach, it becomes difficult to compute multiplicity of submodules. For representations of finite groups, multiplicities are given by some scalar product of characters. Taking inspirations from this theory, the goal of this next section is to explore a similar concept for crossed modules of finite groups.

\section{Characters on representations of a crossed module of finite groups}

\subsection{The structure of the vector space of characters of a crossed module}

Similarly to representations of a finite group, representations of a crossed module of finite groups carry a character theory. The concept was first introduced by P. Bantay in \cite{B}. P. Bantay gives the following definition, which is a direct generalization of the theory of characters for representations of a finite group.

\begin{Definition}[Characters of crossed modules of finite groups \cite{B}]
Let $\xX=(G,X,\mu,\partial)$ be a crossed module of finite groups, and $(V,P,Q)$ an object in $\mM(\xX).$ The character of $(V,P,Q)$ is the $\Cc\mhyphen$valued function
$$
\mapp{\psi}{X\times G}{\Cc}{(m,g)}{\Tr_V(P(m)Q(g)).}
$$
This expression is invariant up to isomorphism of representations of $\xX$.
\end{Definition}

We will write $\Irr(\xX)$ the set of characters of irreducible representations of $\xX.$

P. Bantay also defines what he calls class functions of a crossed module $\xX$. This is a set of $\Cc$-valued functions that follow similar rules to those of characters on a crossed module; characters are class functions. The formal definition is the following.

\begin{Definition}[Class functions of a crossed module \cite{B}]
Let $\xX=(G,X,\mu,\partial)$ be a crossed module of finite groups. A class function of $\xX$ is a $\Cc$-valued function $\psi$ such that
\begin{itemize}
\item For all $x\in X, g\in G,$ $\psi(x,g)=0$ if $g\cdot x\not=x$, 
\item For all $x\in X, g,h\in G, $$\psi(h\cdot x, h g h\mm)=\psi(x,g)$.
\end{itemize}
\end{Definition}

Characters of a crossed module of groups are class functions as, for $g,h\in G, x\in X, (V,P,Q)\in \mM(\xX),$ using the cyclic property of the trace,
\begin{itemize}
\item
\begin{align*} \Tr_V(P(x)Q(g))&=\Tr_V(P(x)^2Q(g))\\&=\Tr_V(P(x)Q(g)P(x))\\&=\Tr_V(P(x)P(g\cdot x)Q(g))\\&=\kr{x,g\cdot x}\Tr_V(P(x)Q(g)).
\end{align*}
\item
\begin{align*}
\Tr_V(P(x)Q(g))&=\Tr_V(Q(h)P(x)Q(g)Q(h\mm))\\&=\Tr_V(P(h\cdot x)Q(hgh\mm).
\end{align*}
\end{itemize}
We will write $\mathscr{C}l(\xX)$ the set of class functions of $\xX$ over the field $\Cc.$

$\mathscr{C}l(\xX)$ is a $\Cc$-vector space of dimension $\sum_{[x]\in X/G}n_x$, where $n_x$ is the number of irreducible characters of the group $\Stab{G}{x}$. It carries the natural scalar product
$$
\forall \psi_1, \psi_2\in \mathscr{C}l(\xX), \langle\psi_1,\psi_2\rangle=\dfrac{1}{\vert G\vert }\dsum_{x\in X, g\in G}\overline{\psi_1(x,g)}\psi_2(x,g).
$$
Here $\overline{z}$ is the complex conjugate of the number $z.$

As a subset of $\textsf{D}(\xX)^*$, $\mathscr{C}l(\xX)$ inherits a richer structure. We define a multiplication that turns $\mathscr{C}l(\xX)$ into a $\Cc$-algebra by the formula
$$
\forall \psi_1, \psi_2\in \mathscr{C}l(\xX), \forall x\in X, g\in G,\quad (\psi_1\psi_2)(x,g)=\dsum_{y\in X}\psi_1(y,g)\psi_2(y\mm x, g).
$$

One can easily check that the multiplication of two class functions is still a class function. For a simple module $I_i^s$ of $\textsf{D}(\xX)$, we associate the character

\begin{align}
\mapp{\psi_{i}^s}{G\times X}{\Cc}{(x,a)}{\Tr_{I_i^s}(m((\kr{x}\ot a)\ot \_))}\label{irrcharmult}
\end{align}

This formula explicitly confirms that irreducible characters are a basis of the vector space of all characters of the crossed module $X$, since $\textsf{D}(\xX)$ is semi-simple and using a property of the trace $\psi_i^s+\psi_j^z=(x,a)\mapsto \Tr_{I_i^s\oplus I_j^z}(m((\kr{x}\ot a)\ot \_))$. Moreover, the (finite) dimension of this vector space is the number of simple $\textsf{D}(\xX)$-modules up to isomorphism (see theorem \ref{nbrsimples}). This is also the dimension of $\mathscr{C}l(\xX)$ and characters are class functions, from this we deduce that the vector space of characters is exactly the vector space of class functions.

From now on, we will write $\psi_i^s$ for the irreducible character that corresponds to some irreducible character $\chi_i$ of $\Stab{G}{s}$. We will also write $\psi_i^s(\kr{x}\ot a):=\psi_i^s(x,a)$, making explicit that characters are elements of $\textsf{D}(\xX)^*$.

Let us write down essential properties of characters, as mentioned in \cite{B}.

\begin{Proposition}
Let $\xX=(G,X,\mu, \partial)$ be a crossed module of finite groups and $\Irr(\xX)$ the set of irreducible characters for its representation theory over $\Cc$.
We have the following properties.
\begin{itemize}
\item $\Irr(\xX)$ is an orthonormal basis of $\mathscr{C}l(\xX)$ for the scalar product $\langle,\rangle$.
\item If $A,B$ are finite dimensional representations of $\xX$, with characters $\chi_A$ and $\chi_B$, then $\chi_A+\chi_B$ is the character of $A\oplus B.$
\item Keeping the same hypotheses, $\chi_A\chi_B$ is the character of $A\ot B$.
\item The multiplicity of $A$ in $B$ is given by $\langle \chi_A, \chi_B\rangle$.
\item The dimension of $A$ is given by $d=\dsum_{x\in X}\chi_A(x, 1)$.
\end{itemize}
\end{Proposition}
\begin{proof}
We prove the first point.

We already saw that $\Irr(\xX)$ was a basis of $\mathscr{C}l(\xX).$ Let us check it is orthonormal. Now, for irreducible characters $\psi_i^s$ and $\psi_j^z$, they can only be non-zero for a single orbit in $X$. For both to be simultaneously non-zero, it is necessary to have $s=z$. Then, $\psi_i^s(g\cdot s,ghg\mm)=\psi_i^s(s,h)$, so
\begin{align*}
\langle \psi_i^s, \psi_j^s\rangle&=\frac{1}{\vert G\vert }\dsum_{g\in G, m\in X}\overline{\psi_i^s(m,g)}\psi_j^s(m,g)\\&=\frac{\vert G\cdot s\vert }{\vert G\vert }\dsum_{g\in \Stab{G}{s}}\overline{\psi_i^s(s,g)}\psi_j^s(s,g)\\&=\frac{1}{\vert \Stab{G}{s}\vert }\dsum_{g\in \Stab{G}{s}}\overline{\psi_i^s(s,g)}\psi_j^s(s,g)\\&=\kr{i,j}.
\end{align*}
The last line comes from \eqref{irrcharmult}.

The next two points come from basic properties of the trace. The fourth point is a direct consequence of the first and second points, and lastly
\begin{align*}
\di{A}&=\Tr_A(\id{A})\\&=\dsum_{m\in X}\Tr_A(P(m)Q(1))\\&=\dsum_{m\in X}\chi_A(m,1).
\end{align*}
\end{proof}

These are all the prerequisites that we need to establish character tables of crossed modules of finite groups. We will now study their general form. Notice that since $\mM(\xX)$ is braided, the multiplication of characters is commutative.
\subsection{Character tables}\label{char}

We define a character table of a crossed module of finite groups as a square table such that the entry $(i,j)$ represents the value of the $i$-th irreducible character evaluated on the $j$-th class of the crossed module $\xX=(G,X, \mu, \partial)$.

We can give a fairly precise description of the character table of a crossed module of finite groups. The tables are diagonal by block tables for which each diagonal block is the character table of the stabilizer of a representative of an orbit.
More visually, we have the following proposition.

\begin{Proposition}
Let $\xX=(G,X,\mu,\partial)$ be a crossed module of finite groups. If $\{s_k\}_{1\leq k\leq n}$ is a set of representatives for orbits of $X$ under the action of $G$, we get the following character table for $\xX$.

\vspace{.5cm}

\begin{center}

\begin{tabular}{|c||c|c|c|c|c|c|
}
\hline
Orbits & \multicolumn{3}{c|}{$s_1$} & $s_2$ & $\ldots$ & $s_n$ \\
\hline\hline
Conjugation classes of $\Cc\Stab{G}{\_}$ & \multicolumn{1}{|c|}{$[1]$}&\multicolumn{2}{|c|}{$\ldots$} &\multicolumn{3}{c|}{$\ldots$} \\
\hline \hline
$\begin{matrix}
\text{Characters of $\xX$} \\
\text{of the form }
\psi_i^{s_1}
\end{matrix}$ & \multicolumn{3}{c|}{$\begin{matrix}
\text{Character table of}\\
\text{the group } \Stab{G}{s_1}
\end{matrix}$}&\multicolumn{3}{c|}{\huge 0}\\

\hline
$\vdots$& \multicolumn{3}{c|}{\huge 0}&\multicolumn{2}{c|}{$\ddots$}&\huge 0\\
\hline
$\begin{matrix}
\text{Characters of $\xX$} \\
\text{of the form }
\psi_i^{s_n}
\end{matrix}$&\multicolumn{5}{c|}{\huge 0} & $\begin{matrix}
\text{Character table of}\\
\text{the group } \Stab{G}{s_n}
\end{matrix}$\\
\hline
\end{tabular}

\vspace{.5cm}

\end{center}

\par

Here $0$ represents a table of the appropriate size, such that each cell is filled with the number $0$.
\end{Proposition}
\begin{proof}
We already proved that an irreducible character is non-zero only for a single orbit of $X$ under the action of $G$, and then only for the stabilizer of elements of this particular orbit for the variable that lives in $G$. Moreover, characters are class functions and, as such, their values are determined only by their values for a set of representatives of orbits.

Now we consider $s\in X$. For $(V, P, Q)$ a finite dimensional irreducible representation of $\xX$ with graduation in $G\cdot s$, there is an induced $\Stab{G}{s}$ representation on $V$ given by $\nu : g\mapsto (v\mapsto Q(g)P(s)v=P(s)Q(g)v)$. This is an irreducible representation of $\Stab{G}{s}$ as if it were not the case, $(V,P,Q)$ would not be irreducible either for $\xX$. Moreover, the character of the irreducible representation $\nu$ is $$\Tr_{P(s)V}(Q(\_))=\Tr_V(Q(\_)P(s))=\Tr_V(P(s)Q(\_))=\psi_V(s,\_)$$ as $P(s)$ is a projection.
\end{proof}

We will give two examples to illustrate, one which is a Drinfeld quantum double and one that is not. Some more example have been computed in \cite{DeghaniDavvaz} and by myself while studying the subject.

\begin{example}[The Drinfeld quantum double of the symmetric group $\Sym_3$]

The conjugacy classes are well known. There are three of them : $[1], [(12)]$ and $[(123)]$.
Their respective associated stabilizer, which here correspond to the centralizer of a representative, are $\Sym_3$, $\{1, (12)\}\cong \Zz/2\Zz$ and $\{1, (123), (132)\}\cong \Zz/3\Zz$. Their respective table of characters are also well known, and therefore the table of characters of $\xX=(\Sym_3, \Sym_3, \conj{\_},\id{\Sym_3})$ is the following. Here $\omega$ is a primitive third root of unity in $\Cc$. The simple representations of this Drinfeld double have been fully described in \cite{L}.

\centering

\vspace{.5cm}
\begin{tabular}{|c||c|c|c!{\vrule width 2pt}c|c!{\vrule width 2pt}c|c|c|}
\hline
Orbit & \multicolumn{3}{c!{\vrule width 2pt}}{$\id{}$} & \multicolumn{2}{c!{\vrule width 2pt}}{$(12)$} &\multicolumn{3}{c|}{$(123)$}\\
\hline \hline
Conjugacy class & $[1]$ &$[(12)]$&$[(123)]$&$[1]$&$[(12)]$&$[1]$&$[(123)]$&$[(132)]$\\
\hline \hline
$\psi_1^1$& 1&1&1&0&0&0&0&0 \\
\hline
$\psi_2^1$&1&-1&1&0&0&0&0&0\\
\hline
$\psi_3^1$&2&0&-1&0&0&0&0&0\\
\Xhline{2pt}
$\psi_1^{(12)}$&0&0&0&1&1&0&0&0\\
\hline
$\psi_2^{(12)}$&0&0&0&1&-1&0&0&0\\
\Xhline{2pt}
$\psi_1^{(123)}$&0&0&0&0&0&1&1&1\\
\hline
$\psi_2^{(123)}$&0&0&0&0&0&1&$\omega$&$\omega^2$\\
\hline
$\psi_3^{(123)}$&0&0&0&0&0&1&$\omega^2$&$\omega$\\
\hline
\end{tabular}
\vspace{.5cm}

\end{example}

\begin{example}
We look at $\xX=(\Aut{D_4}\cong D_4, D_4, \text{ev}, x\mapsto x\_x\mm)$, where $D_4=\langle r,s\vert s^2=1, r^4=1, r^3s=sr \rangle$ is the dihedral group of symmetries of the square. This is not a Drinfeld quantum double even though $\Aut{D_4}\cong D_4$. We only need to specify that $\Stab{\Aut{D_4}}{s}\cong \Zz/2\Zz$, $\Stab{\Aut{D_4}}{r}\cong \Zz/ 4\Zz$ and $\Stab{\Aut{D_4}}{r^2}\cong\Stab{\Aut{D_4}}{1}\cong D_4$. The character table of $\xX$ is therefore the following.

\vspace{.5cm}
\centering

\begin{tabular}{|c||c|c|c|c|c!{\vrule width 2pt}c|c|c|c|c!{\vrule width 2pt}c|c|c|c!{\vrule width 2pt}c|c|}
\hline
Orbit & \multicolumn{5}{c!{\vrule width 2pt}}{$1$} & \multicolumn{5}{c!{\vrule width 2pt}}{$r^2$}&\multicolumn{4}{c!{\vrule width 2pt}}{$r$}& \multicolumn{2}{c|}{$s$}\\
\hline
\hline
Conjugacy class & $[1]$& $[r]$& $[s]$& $[r^2]$& $[sr]$&$[1]$& $[r]$& $[s]$& $[r^2]$& $[sr]$&$[1]$& $[r]$& $[r^2]$& $[r^3]$& $[1]$& $[s]$\\
\hline
\hline
$\psi_1^1$& 1& 1& 1& 1& 1& 0& 0& 0& 0& 0& 0& 0& 0& 0& 0&0\\
\hline
$\psi_2^1$& 1& -1& 1& 1& -1& 0& 0& 0& 0& 0& 0& 0& 0& 0& 0&0\\
\hline
$\psi_3^1$& 1& -1& -1& 1& 1& 0& 0& 0& 0& 0& 0& 0& 0& 0& 0&0\\
\hline
$\psi_4^1$& 1& 1& -1& 1& -1& 0& 0& 0& 0& 0& 0& 0& 0& 0& 0&0\\
\hline
$\psi_5^1$& 2& 0& 0& -2& 0& 0& 0& 0& 0& 0& 0& 0& 0& 0& 0&0\\
\Xhline{2pt}
$\psi_1^{r^2}$& 0& 0& 0& 0& 0& 1& 1& 1& 1& 1& 0& 0& 0& 0& 0&0\\
\hline
$\psi_2^{r^2}$& 0& 0& 0& 0& 0& 1& -1& 1& 1& -1& 0& 0& 0& 0& 0&0\\
\hline
$\psi_3^{r^2}$& 0& 0& 0& 0& 0& 1& -1& -1& 1& 1& 0& 0& 0& 0& 0&0\\
\hline
$\psi_4^{r^2}$& 0& 0& 0& 0& 0& 1& 1& -1& 1& -1& 0& 0& 0& 0& 0&0\\
\hline
$\psi_5^{r^2}$& 0& 0& 0& 0& 0& 2& 0& 0& -2& 0& 0& 0& 0& 0& 0&0\\
\Xhline{2pt}
$\psi_1^r$& 0& 0& 0& 0& 0& 0& 0& 0& 0& 0& 1& 1& 1& 1& 0&0\\
\hline
$\psi_2^r$& 0& 0& 0& 0& 0& 0& 0& 0& 0& 0& 1& i& -1& -i& 0&0\\
\hline
$\psi_3^r$& 0& 0& 0& 0& 0& 0& 0& 0& 0& 0& 1& -1& 1& -1& 0&0\\
\hline
$\psi_4^r$& 0& 0& 0& 0& 0& 0& 0& 0& 0& 0& 1& -i& -1& i& 0&0\\
\Xhline{2pt}
$\psi_1^s$& 0& 0& 0& 0& 0& 0& 0& 0& 0& 0& 0& 0& 0& 0& 1&1\\
\hline
$\psi_2^s$& 0& 0& 0& 0& 0& 0& 0& 0& 0& 0& 0& 0& 0& 0& 1&-1\\
\hline

\end{tabular}
\par
\end{example}

\vspace{1cm}

\subsection{The Clebsch-Gordan formula via character theory}\label{TheClebschGordanFormulaViaCharacterTheory}

We will begin by reminding the following theorem, from P. Bantay.

\begin{theorem}[P. Bantay, \cite{B}]
For a crossed module of finite groups $\xX=(G,X,\mu,\partial)$, the multiplicity of a finite dimensional simple $\textsf{D}(\xX)$-module $I_k^t$ in $I_i^s\ot I_j^z$ is
\begin{align*}
N_{(s,i),(z,j)}^{(t,k)}=\langle \psi_i^s\psi_j^z,\psi_k^t\rangle.
\end{align*}
\end{theorem}

We can improve this formula with the previous considerations. We use the fact that an irreducible character is non-zero only for a single orbit.

\begin{align*}
N_{(s,i),(z,j)}^{(t,k)}&=\frac{1}{|G|}\dsum_{g\in G, n,m\in X}\psi_i^s(n,g)\psi_j^z(n\mm m,g)\overline{\psi_k^t(m,g)}\\
&=\frac{1}{\vert G\vert }\dsum_{g\in G,n,m\in X}\psi_i^s(n,g)\psi_j^z(m,g)\overline{\psi_k^t(nm,g)}\\
&=\dfrac{1}{\vert G\vert \vert \Stab{G}{s}\vert \vert \Stab{G}{z}\vert }\dsum_{a,b,g\in G}\psi_i^s(a\cdot s,g)\psi_j^z(b\cdot z,g)\overline{\psi_k^t((a\cdot s)(b\cdot z), g)}\\
&=\frac{1}{\vert G\vert \vert \Stab{G}{s}\vert \vert \Stab{G}{z}\vert }\dsum_{a,b,g\in G}\psi_i^s(s,a\mm ga)\psi_j^z(z,b\mm gb)\overline{\psi_k^t(s(a\mm b)\cdot z, a\mm ga)}
\end{align*}
\begin{align*}
&=\frac{1}{\vert G\vert \vert \Stab{G}{s}\vert \vert \Stab{G}{z}\vert }\dsum_{a,b,g\in G}\psi_i^s(s,g)\psi_j^z(z,b\mm aga\mm b)\overline{\psi_k^t(s(a\mm b)\cdot z, g)}\\
&=\frac{1}{\vert \Stab{G}{s}\vert \vert \Stab{G}{z}\vert }\dsum_{a,g\in G}\psi_i^s(s,g)\psi_j^z(z,a\mm ga)\overline{\psi_k^t(s(a\cdot z),g)}
\end{align*}

From as soon as the second line, the condition we have already found about $(G\cdot s)(G\cdot z)\cap (G\cdot t)\not=\emptyset$ appears, as elements in the sum are non-zeros only if the variables $n,m$ are elements of $G\cdot s$, $G\cdot z$ respectively and their product is in $G\cdot t$. We can only notice that the last sum is only over $G$ and does not directly involve $X$, in a way reminiscent of what we proved in \ref{Cleb-Gor}. This result is limited in terms of understanding. Previously we were able to interpret beyond equations what it meant for a simple module to be a submodule of the tensor product of two simple modules, here it seems complicated to get a satisfying interpretation. However, this formula could be implemented in a computer, and should be more efficient for effective computations.

These coefficients appear, as stated in \cite{DVVV} for the case of a Drinfeld double, as fusion coefficients for the Grothendieck algebra $K_X(G)$ of $G\mhyphen$equivariant complex vector bundles over the finite group $X$. This remark ties our category to orbifold conformal field theories, $(K_X(G), \oplus, \ot )$ becoming the fusion algebra.

Consider $I_i^s, I_j^z$ two complex finite dimensional irreducible representations of $\xX=(G,X, \mu, \partial )$ a crossed module of finite groups. Both are elements of $K_X(G)$, and, paraphrasing relations above,
\begin{align*}
I_i^s\ot I_j^z=\bigoplus_{t\in \hH}\bigoplus_{1\leq k \leq n_t} (I_k^t)^{\oplus N_{(s,i),(z,j)}^{(t,k)}}.
\end{align*}

\section{Quantum invariants}\label{LastSect}

We proved in previous sections that the category $\mM(\xX)$ has a lot of structure. It is a category of finite-dimensional modules over a braided finite dimensional semi-simple Hopf algebra, which makes it a fusion category.

Moreover, the dual of an object $M$, with basis $\bB$, in $\mM(\xX)$, is the dual vector space $M^*$ together with the left $\textsf{D}(\xX)$-module structure on a basis
$$\mapp{\phi}{\textsf{D}(\xX)\ot M^*}{M^*}{\kr{x}\ot g\ot \kr{m}}{\dsum_{n \in \bB}\kr{m}((\kr{g\mm \cdot x\mm }\ot g\mm)\cdot n )\kr{n}.}$$

The evaluation and coevaluation morphisms are the usual ones, and turn $\mM(\xX)$ into a rigid category.

Lastly, we give the following theorem. We refer to \cite{K}, Chapter XIV for definitions of ribbon element or ribbon algebra.

\begin{theorem}
For a crossed module of finite groups $\xX=(G,X,\mu,\partial)$, $\textsf{D}(\xX)$ is a ribbon algebra with ribbon element
$$
\theta=\dsum_{x\in X}\kr{x}\ot \partial(x\mm).
$$
\end{theorem}
\begin{proof}
The ribbon element is actually exactly the Drinfeld element of $\textsf{D}(\xX)$ (\cite{Rad}, Definition 12.2.9). The fact that $\theta$ is central comes from the property of the Drinfeld element (\cite{Rad}, VIII.4), that for all $y$ in $\textsf{D}(\xX)$, $S\circ S(y)=\theta y \theta\mm$, and $S\circ S=\id{\textsf{D}(\xX)}$. $\varepsilon(\theta)=1$ can be directly computed or obtained as a consequence of $\theta$ being the Drinfeld element of $\textsf{D}(\xX)$, the equation $\Delta(\theta)=(\theta\ot\theta)(R_{21}R)\mm$ is also always satisfied for the Drinfeld element, see for example proposition 4.3 in \cite{Ohtsuki}. Through direct computations we also get 
\begin{align*}
    S(\theta)&=\sum_{x\in X}\kr{\partial(x)\cdot x\mm}\ot \partial(x)\\&=\sum_{x\in X}\kr{x\mm}\ot \partial(x)\\&=\theta.
\end{align*}
This obviously implies $\theta^2=S(\theta)\theta$, and once again since $\theta$ is the Drinfeld element of $\textsf{D}(\xX)$, we conclude that it is indeed also a ribbon element of $\textsf{D}(\xX)$.
\end{proof}

From $\theta$ we can deduce all other possible ribbon element (\cite{Rad}, Theorem 12.3.6.). Follows a corollary useful for generating less trivial invariants later.

\begin{Corollary}
Ribbon elements of $\textsf{D}(\xX)$ are indexed by the subgroup $C$ of $G$ of central elements of order 2 that act trivially on $X$. The set of all ribbon elements is $\{(1 \ot c)\theta | c \in C\}$.
\end{Corollary}

Because $\textsf{D}(\xX)$ is a ribbon algebra, $\mM(\xX)$ is a ribbon category, and using Rechetikhin-Turaev's theorem (\cite{Turaev}, Chapter I Theorem 2.5.), we are able to generate a ribbon invariant for any $c\in C$. By definition, this invariant is a monoidal functor $F$ from $\textsf{Rib}_{\mM(\xX)}$ the category of $\mM(\xX)\mhyphen$coloured ribbon tangles (as defined in \cite{Turaev}, Chapter I. Subsection 2.2.) to $\mM(\xX)$. For elementary tangles, the functor has the following values.  The following theorem is an instance of the aforementioned Rechetikhin-Turaev\rq{}s theorem, together with our previous study of the Hopf algebra $\textsf{D}(\xX)$.

\begin{theorem}
For a crossed module of finite groups $\xX=(G,X,\mu,\partial)$, the Hopf algebra $\textsf{D}(\xX)$ of \ref{DefDX} generates a ribbon invariant $F_c$, for any $c\in C$, whose values over elementary tangles are defined as follows. 
Let $(V, P_V, Q_V), (W, P_W, Q_W) $ be objects in $\mM(\xX)$. The blue strand is $V$-coloured and the red one is $W$-coloured. Morphisms are written in a basis $\bB_V$ (resp. $\bB_W)$ of the vector space $V$ (resp. $W$) (or the corresponding dual basis).

\begin{tikzpicture}[arrow/.style={->}, rarrow/.style={<-}]
\begin{knot}[consider self intersections, end tolerance=3pt]
\strand[blue, rarrow] (0,0) to[out=up,in=down] ++(1,1);
\strand[red, arrow] (0,1) to[out=down,in=up] ++(1,-1);
\strand[red, arrow] (0,-1) to[out=down,in=up] ++(1,-1);
\strand[blue, rarrow] (0,-2) to[out=up,in=down] ++(1,1);
\strand[blue, rarrow] (0,-4) to[out=up,in=down] ++(1,1);
\strand[red, rarrow] (0,-3) to[out=down,in=up] ++(1,-1);
\strand[red, rarrow] (0,-5) to[out=down,in=up] ++(1,-1);
\strand[blue, rarrow] (0,-6) to[out=up,in=down] ++(1,1);
\strand[blue, arrow] (0,-7) to[out=down,in=left](0.5,-8) to[out=right, in=down] (1,-7);
\strand[blue, arrow] (0,-10) to[out=up,in=left](0.5,-9) to[out=right, in=up] (1,-10);
\strand[blue, arrow] (1,-11)
to[out=up, in=down] (1,-11.7)
to[out=down, in=right] (0.5,-12.5)
to[out=left, in=down] (0,-12)
to[out=up, in=left] (0.5,-11.5)
to[out=right, in=up] (1,-12.3)
to[out=up, in=up] (1,-13);
\node at (4,0.5) {$\overset{F_c}{\longmapsto}$}; \node at (10,0.5) {$\mapp{c_{V,W}}{V\ot W}{W\ot V}{v\ot w}{\dsum_{x\in X}Q_W(\partial(x))w\ot P_V(x)v}$};
\node at (4,-1.5) {$\overset{F_c}{\longmapsto}$}; \node at (10,-1.5) {$\mapp{c_{W,V}^{-1}}{V\ot W}{W\ot V}{v\ot w}{\dsum_{x\in X}Q_W(\partial(x))\mm w\ot P_V(x)v}$};
\node at (4,-3.5) {$\overset{F_c}{\longmapsto}$}; \node at (10,-3.5) {$\mapp{c_{V,W^*}}{V\ot W^*}{W^*\ot V}{v\ot \kr{w}}{\dsum_{x\in X}\kr{w}(Q_W(\partial(x))\cdot\_) \ot P_V(x)v}$};
\node at (4,-5.5) {$\overset{F_c}{\longmapsto}$}; \node at (10,-5.5) {$\mapp{c_{W^*,V}^{-1}}{V\ot W^*}{W^*\ot V}{v\ot \kr{w}}{\dsum_{x\in X}\kr{w}(Q_W(\partial(x))\mm\cdot\_) w\ot P_V(x)v}$};
\node at (4,-7.5) {$\overset{F_c}{\longmapsto}$}; \node at (10,-7.5) {$\mapp{\text{coev}_V}{\Cc}{V^*\ot V}{1}{\dsum_{v\in \bB_V}\kr{v}\ot v}$};
\node at (4,-9.5) {$\overset{F_c}{\longmapsto}$}; \node at (10,-9.5) {$\mapp{\text{ev}_V}{V^*\ot V}{\Cc}{\kr{v}\ot v'}{\kr{v,v'}}$};
\node at (4,-11.5) {$\overset{F_c}{\longmapsto}$}; \node at (10,-11.5) {$\mapp{\theta_V}{V}{V}{v}{\dsum_{x\in X}P_V(x)Q_V(\partial(x))\mm Q_V(c)v}$};
\end{knot}
\end{tikzpicture}

These values are enough to derive any value of $F_c$ for any ribbon of $\textsf{Rib}_{\mM(\xX)}$ (\cite{Turaev}, Chapter I. Lemma 3.1.1.). 
\end{theorem}

Since the category $\mM(\xX)$ is not modular, unless $\textsf{D}(\xX)$ is the Drinfeld double of $G$ (see \cite{MaierSchweigert}), we cannot as easily access stronger invariants such as those discussed in Chapter III of \cite{Turaev} without loss of generality.

\vspace{.5cm}
\textbf{Declarations : }

\textbf{Ethical Approval } Not applicable

\textbf{Funding } This work was financed by Laboratoire de mathématiques Jean Leray.

\textbf{Availability of data and materials } Not applicable.

\bibliographystyle{plain}
\bibliography{Biblio}

Ony Aubril\newline
Université Catholique de Louvain-la-Neuve\newline
1348 Louvain-la-Neuve\newline
Belgium \newline
e-mail : ony.aubril@uclouvain.be
\end{document}